\documentclass{CSML}

\def\dOi{12(4:3)2016}
\lmcsheading%
{\dOi}
{1--11}
{}
{}
{Sep.~17, 2015}
{Oct.~12, 2016}
{}

\ACMCCS{[{\bf Theory of computation}]: Models of computation; Design
  and analysis of algorithms [{\bf Mathematics of computing}]:
  Mathematical analysis---Numerical analysis; Continuous
  mathematics---Continuous functions}


\usepackage{hyperref}
\usepackage{hyphenat}
\usepackage{amssymb}
\theoremstyle{plain}
\theoremstyle{plain}
\theoremstyle{plain}\newtheorem{example}[thm]{Example}
\theoremstyle{plain}\newtheorem{remark}[thm]{Remark}
\theoremstyle{plain}
\begin{document}

\title[FROM LOCAL TO GLOBAL PRESENCE OF PROPERTIES OF FUNCTIONS]{SOME THEOREMS ON PASSING\\FROM LOCAL TO GLOBAL PRESENCE\\OF PROPERTIES OF FUNCTIONS}

\author[D.~Skordev]{Dimiter Skordev}	%required
\address{Sofia University, Faculty of Mathematics and Informatics, Sofia, Bulgaria}	%required
\email{skordev@fmi.uni-sofia.bg}  %optional
%\thanks{}	%optional

%% etc.

%% required for running head on odd and even pages, use suitable
%% abbreviations in case of long titles and many authors:

%% mandatory lists of keywords and classifications:
%\titlecomment{}
%%%%%%%%%%%%%%%%%%%%%%%%%%%%%%%%%%%%%%%%%%%%%%%%%%%%%%%%%%%%%%%%%%%%%%%%%%%

%% the abstract has to PRECEED the command \maketitle:
%% be sure not to issue the \maketitle command twice!

\begin{abstract}
When given a class of functions and a finite collection of sets, one might be interested whether the class in question contains any function whose domain is a subset of the union of the sets of the given collection and whose restrictions to all of them belong to this class. The collections with the formulated property are said to be {\em strongly join permitting} for the given class (the notion of join permitting collection is defined in the same way, but without the words ``a subset of''). Three theorems concerning certain instances of the problem are proved. A necessary and sufficient condition for being strongly join permitting is given for the case when, for some $n$, the class consists of the potentially partial recursive functions of $n$ variables, and the collection consists of sets of $n${\hyp}tuples of natural numbers. The second theorem gives a sufficient condition for the case when the class consists of the continuous partial functions between two given topological spaces, and the collection consists of subsets of the first of them (the condition is also necessary under a weak assumption on the second one). The third theorem is of a similar character but, instead of continuity, it concerns computability in the spirit of the one in effective topological spaces.
\end{abstract}

\maketitle

%% start the paper here:
\section{Introduction}\label{S:intro}
To show that a given function $f$ has a certain property, it is sometimes useful to cover the domain of $f$ with an appropriate collection of sets and to show that the property in question is present for the restriction of $f$ to any of the sets from this collection. For instance, suppose $f$ is the $\arctan$ function considered for all real values of its argument, and we aim at proving the computability of this function (in the sense of the computable analysis). To do this, we may firstly consider the restrictions of $f$ to the intervals $[-1,1]$, $[1,+\infty)$ and $(-\infty,-1]$. These restrictions can be shown to be computable as follows. The function $f\!\upharpoonright[-1,1]$ is computable thanks to the fact that
\begin{equation}\label{ml}
\arctan x=x-\frac{x^3}{3}+\frac{x^5}{5}-\frac{x^7}{7}+\cdots
\end{equation}
for all real numbers $x$ with $|x|\le 1$. The computability of $f\!\upharpoonright[1,+\infty)$ can be derived from here by using that
$$\arctan x=\frac{\pi}{2}-\arctan\frac{1}{x}$$
for all positive real numbers $x$, and the computability of $f\!\upharpoonright(-\infty,-1]$ can be reduced to the one of $f\!\upharpoonright[1,+\infty)$ by means of the equality
$$\arctan x=-\arctan(-x).$$
On the other hand, two successive applications of Lemma 4.3.5 from \cite{ca} (after correcting the misprints in the lemma by replacement of $a$ with $c$) show the computability of any real function whose restrictions to the above-mentioned three intervals are computable.

\begin{defi}\label{jp}
Let $\mathbf{C}$ be a class of functions. A collection $\mathcal{A}$ of sets will be said to be {\em join permitting for $\mathbf{C}$} if the following implication holds for any function $f$ whose domain is~$\bigcup\mathcal{A}$:
\begin{equation}\label{jpr}
\forall A\in\mathcal{A}(f\!\upharpoonright_A\in\mathbf{C})\Longrightarrow f\in\mathbf{C};
\end{equation}
the collection $\mathcal{A}$ will be said to be {\em strongly join permitting for $\mathbf{C}$} if the implication \eqref{jpr} holds for any function $f$ whose domain is a subset of $\bigcup\mathcal{A}$. 
\end{defi}

We will be interested in the above notions in the case of finite collections $\mathcal{A}$, i.e.\@ in the case when $\mathcal{A}$ is a finite set of sets.

\begin{example}\label{intervals}
{\em If $\mathbf{C}$ is the class of all computable partial functions from $\mathbb{R}$ to $\mathbb{R}$, and $c_1<c_2<\cdots<c_{r-1}<c_r$ are computable real numbers then the collection of the intervals $(-\infty,c_1],[c_1,c_2],\ldots,[c_{r-1},c_r],[c_r,+\infty)$ is strongly join permitting for $\mathbf{C}$ (this can be shown by means of $r$ successive applications of the above-mentioned lemma from \cite{ca}).}
\end{example}

By trivial reasons, if some of the sets in a given collection of sets is the domain of no function of the class $\mathbf{C}$ then the collection turns out to be join permitting for $\mathbf{C}$, and if some of these sets contains the domain of no function of $\mathbf{C}$ then the collection is strongly join permitting for $\mathbf{C}$.

There are simple examples of classes $\mathbf{C}$ such that finite collections exist which are join permitting for $\mathbf{C}$ without being strongly join permitting for $\mathbf{C}$.\footnote{For instance, let $A_1$ and $A_2$ be sets such that $A_1\setminus A_2$, $A_2\setminus A_1$ and $A_1\cap A_2$ are non-empty, and let $\mathbf{C}$ consist of the partial functions in $A_1\cup A_2$ with one-element ranges. Then the collection $\{A_1,A_2\}$ is join permitting for $\mathbf{C}$, but it is not strongly join permitting for $\mathbf{C}$.} We will direct our attention to the notion of strongly join permitting collection.

For any set $A$, let $A_\mathbf{C}$ be the set of the elements of $A$ which belong to the domain of at least one function of $\mathbf{C}$. The set $A$ will be called {\em relevant for $\mathbf{C}$} if $A_\mathbf{C}=A$. Clearly, $A_\mathbf{C}$ is always relevant for $\mathbf{C}$.

\begin{lem}\label{relev}
Let $\mathbf{C}$ be a class of functions, and $\mathcal{A}$ be a collection of sets. The collection $\mathcal{A}$ is strongly join permitting for $\mathbf{C}$ iff the corresponding collection $\{A_\mathbf{C}\,|\,A\in\mathcal{A}\}$ is strongly join permitting for $\mathbf{C}$.
\end{lem}

\proof 
Let us firstly suppose that $\mathcal{A}$ is strongly join permitting for $\mathbf{C}$. Let $f$ be a function such that $\mathrm{dom}(f)\subseteq\bigcup\{A_\mathbf{C}\,|\,A\in\mathcal{A}\}$ and $f\!\!\upharpoonright_{A_\mathbf{C}}\in\mathbf{C}$ for any $A\in\mathcal{A}$. We will show that $f\in\mathbf{C}$. The domain of $f$ is a subset of $\bigcup\mathcal{A}$ because $A_\mathbf{C}\subseteq A$ for any set $A$. Any point of $\mathrm{dom}(f)$ belongs to $A_\mathbf{C}$ for some $A\in\mathcal{A}$, hence it belongs to the domain of some function from~$\mathbf{C}$. Thus $f\!\!\upharpoonright_A=f\!\!\upharpoonright_{A_\mathbf{C}}$ and consequently $f\!\!\upharpoonright_A\in\mathbf{C}$ for any $A\in\mathcal{A}$. By the assumption that $\mathcal{A}$ is strongly join permitting for $\mathbf{C}$, this implies $f\in\mathbf{C}$.

For the reasoning in the opposite direction, suppose the collection $\{A_\mathbf{C}\,|\,A\in\mathcal{A}\}$ is strongly join permitting for $\mathbf{C}$. Let $f$ be a function such that $\mathrm{dom}(f)\subseteq\bigcup\mathcal{A}$ and $f\!\!\upharpoonright_A\in\mathbf{C}$ for any $A\in\mathcal{A}$. We will show that $f\in\mathbf{C}$. For any $A\in\mathcal{A}$, the points of $A\cap\mathrm{dom}(f)$ belong to the domain of the function $f\!\!\upharpoonright_A$, which is a function from~$\mathbf{C}$, hence they belong to the subset $A_\mathbf{C}$ of~$A$. Therefore $f\!\!\upharpoonright_A=f\!\!\upharpoonright_{A_\mathbf{C}}$ and consequently $f\!\!\upharpoonright_{A_\mathbf{C}}\in\mathbf{C}$ for any $A\in\mathcal{A}$. Since any point of~$\mathrm{dom}(f)$ belongs to some set from $\mathcal{A}$, we see also that $\mathrm{dom}(f)\subseteq\bigcup\{A_\mathbf{C}\,|\,A\in\mathcal{A}\}$. By the assumption that $\{A_\mathbf{C}\,|\,A\in\mathcal{A}\}$ is strongly join permitting for $\mathbf{C}$, this implies $f\in\mathbf{C}$.
\qed

For some classes $\mathbf{C}$, a trivial characterization is possible of the finite collections which are strongly join permitting for $\mathbf{C}$. For instance, such is the case when $\mathbf{C}$ is the class of the one-argument partial recursive functions -- any finite collection of sets is strongly join permitting for $\mathbf{C}$ in this case (because the class in question consists of the functions whose graphs are recursively enumerable subsets of~$\mathbb{N}^2$, and the union of finitely many recursively enumerable subsets of $\mathbb{N}^2$ is also recursively enumerable). However, the problem is more complicated for certain other classes. For some such classes $\mathbf{C}$ and the finite collections of sets relevant for them, we will give necessary and sufficient conditions for being strongly join permitting for $\mathbf{C}$, and these conditions will be in a similar spirit. The following definition will be used:

\iffalse\begin{defi}\label{cup}
Whenever $\mathcal{A}$ is a collection of sets, and $\mathcal{K}$ is a subcollection of $\mathcal{A}$, the union $\bigcup_{A\in\mathcal{K}}A$ will be denoted by $\bigcup\mathcal{K}$.
\end{defi}\fi

\begin{defi}\label{dsep}
A set $H$ will be said to {\em separate} a set $P$ from a set $Q$ if $H\supseteq P$ and $H\cap Q=\varnothing$. 
\end{defi}

Of course, a set separating $P$ from $Q$ exists iff $P\cap Q=\varnothing$, but this equivalence can turn out to be no more valid if some restrictions on the separating set are imposed.

\begin{lem}\label{lss}
Let $\mathcal{A}$ be a collection of sets, $\mathcal{K}$ be a subcollection of $\mathcal{A}$ and $H$ be a set which separates $\bigcup\mathcal{K}\setminus\bigcup(\mathcal{A}\setminus\mathcal{K})$ from $\bigcup(\mathcal{A}\setminus\mathcal{K})\setminus\bigcup\mathcal{K}$. Then:
\begin{enumerate}[label=\({\alph*}]
\item $(\bigcup\mathcal{A})\cap H\subseteq\bigcup\mathcal{K}$; 
\item if $x\in\bigcup\mathcal{A}$ and $\mathcal{K}\supseteq\{A\in\mathcal{A}\,|\,x\in A\}$ then $x\in H$.
\end{enumerate}
\end{lem}

\proof
By the assumptions of the lemma, $H\supseteq\bigcup\mathcal{K}\setminus\bigcup(\mathcal{A}\setminus\mathcal{K})$ and \mbox{$H\cap\left(\bigcup(\mathcal{A}\setminus\mathcal{K})\setminus\bigcup\mathcal{K}\right)=\varnothing$.} To prove (a), suppose $x\in\left(\bigcup\mathcal{A}\right)\cap H$. Then $x\in A$ for some $A\in\mathcal{A}$. If $A\in\mathcal{K}$ then $x\in\bigcup\mathcal{K}$. In the opposite case, $x\in\bigcup(\mathcal{A}\setminus\mathcal{K})$ and then surely $x\in\bigcup\mathcal{K}$ again, since otherwise $x$ would belong to $\bigcup(\mathcal{A}\setminus\mathcal{K})\setminus\bigcup\mathcal{K}$, and this is impossible since $x\in H$. To prove~(b), suppose $x\in\bigcup\mathcal{A}$ and $\mathcal{K}\supseteq\{A\in\mathcal{A}\,|\,x\in A\}$. Then clearly $x\in\bigcup\mathcal{K}$. On the other hand, $x\not\in\bigcup(\mathcal{A}\setminus\mathcal{K})$, because the opposite would imply that $x\in A$ for some $A\in\mathcal{A}$ such that $x\not\in A$. Thus $x\in\bigcup\mathcal{K}\setminus\bigcup(\mathcal{A}\setminus\mathcal{K})$ and therefore $x\in H$.
\qed

\begin{lem}\label{const}
Let $\mathcal{A}$ be a collection of sets, $\mathcal{K}$ be a subcollection of $\mathcal{A}$, and $c_1,c_2$ be two given objects. Let $f:\left(\bigcup\mathcal{K}\setminus\bigcup(\mathcal{A}\setminus\mathcal{K})\right)\cup\left(\bigcup(\mathcal{A}\setminus\mathcal{K})\setminus\bigcup\mathcal{K}\right)\to\{c_1,c_2\}$ be such that $f(x)=c_1$ for any $x\in\bigcup\mathcal{K}\setminus\bigcup(\mathcal{A}\setminus\mathcal{K})$ and $f(x)=c_2$ for any $x\in\bigcup(\mathcal{A}\setminus\mathcal{K})\setminus\bigcup\mathcal{K}$. Then $\mathrm{dom}(f)\subseteq\bigcup\mathcal{A}$ and, for any $A\in\mathcal{A}$, $f\!\!\upharpoonright_A$ is a constant function.
\end{lem}

\proof
The inclusion $\mathrm{dom}(f)\subseteq\bigcup\mathcal{A}$ is obvious. To prove the other statement of the lemma, suppose $A\in\mathcal{A}$. If $A\in\mathcal{K}$ then $A\subseteq\bigcup\mathcal{K}$ , hence the intersection of $A$ with $\bigcup(\mathcal{A}\setminus\mathcal{K})\setminus\bigcup\mathcal{K}$ is empty and therefore $A\cap\mathrm{dom}(f)\subseteq\bigcup\mathcal{K}\setminus\bigcup(\mathcal{A}\setminus\mathcal{K})$. Similarly, if $A\in\mathcal{A}\setminus\mathcal{K}$ then $A\subseteq\bigcup(\mathcal{A}\setminus\mathcal{K})$, hence the intersection of $A$ with $\bigcup\mathcal{K}\setminus\bigcup(\mathcal{A}\setminus\mathcal{K})$ is empty and therefore $A\cap\mathrm{dom}(f)\subseteq\bigcup(\mathcal{A}\setminus\mathcal{K})\setminus\bigcup\mathcal{K}$.
\qed

\section{Strongly join permitting collections for the class \\of the potentially partial recursive $n${\hyp}ary functions}\label{S:ppr}
Let $n$ be a positive integer, and $\mathbf{C}$ be the class of all potentially partial recursive $n${\hyp}ary functions (as usual, a partial function from $\mathbb{N}^n$ to $\mathbb{N}$ will be called {\em potentially partial recursive} if it is a restriction of some $n${\hyp}ary partial recursive function). 

\begin{thm}\label{ppr}
Let $\mathcal{A}$ be a finite collection of subsets of $\mathbb{N}^n$. The collection $\mathcal{A}$ is strongly join permitting for~$\mathbf{C}$ iff, for each subcollection $\mathcal{K}$ of $\mathcal{A}$, some recursively enumerable subset $H_\mathcal{K}$ of~$\mathbb{N}^n$ separates $\bigcup\mathcal{K}\setminus\bigcup(\mathcal{A}\setminus\mathcal{K})$ from $\bigcup(\mathcal{A}\setminus\mathcal{K})\setminus\bigcup\mathcal{K}$.
\end{thm}

\proof Suppose that, for any subcollection $\mathcal{K}$ of $\mathcal{A}$, some recursively enumerable subset $H_\mathcal{K}$ of~$\mathbb{N}^n$ separates $\bigcup\mathcal{K}\setminus\bigcup(\mathcal{A}\setminus\mathcal{K})$ from $\bigcup(\mathcal{A}\setminus\mathcal{K})\setminus\bigcup\mathcal{K}$. We will prove that $\mathcal{A}$ is strongly join permitting for $\mathbf{C}$. Let $f$ be a function such that $\mathrm{dom}(f)\subseteq\bigcup\mathcal{A}$ and $f\!\upharpoonright_A\in\mathbf{C}$ for any $A\in\mathcal{A}$. Clearly, all values of $f$ are natural numbers. For any $A\in\mathcal{A}$, the function $f\!\upharpoonright_A$ is a restriction of some $n${\hyp}ary partial recursive function $\varphi_A$. Consider an arbitrary element $x$ of $\mathrm{dom}(f)$ and an arbitrary natural number~$y$. Then $x\in\bigcup\mathcal{A}$. We will show that the equality $f(x)=y$ is equivalent to the following condition:
\begin{equation}\label{l}
\exists\mathcal{K}\subseteq\mathcal{A}(x\in H_\mathcal{K}\ \&\ \forall A\in\mathcal{K}(\varphi_A(x)=y)).
\end{equation}
Firstly, suppose that $x\in H_\mathcal{K}\ \&\ \forall A\in\mathcal{K}(\varphi_A(x)=y)$ for some subcollection $\mathcal{K}$ of $\mathcal{A}$. The statement~(a) of Lemma \ref{lss} implies that $x\in A$ for some $A\in\mathcal{K}$ and therefore \mbox{$f(x)=\varphi_A(x)=y$ for this~$A$.} For the reasoning in the opposite direction, suppose that $f(x)=y$. The statement~(b) of Lemma \ref{lss} implies that $x\in H_\mathcal{K}$, where \mbox{$\mathcal{K}=\{A\in\mathcal{A}\,|\,x\in A\}$.} Clearly $\varphi_A(x)=f(x)=y$ for any $A$ belonging to this $\mathcal{K}$.

Thanks to the finiteness of $\mathcal{A}$, the condition (\ref{l}) defines a recursively enumerable subset of~$\mathbb{N}^{n+1}$. An application of the Uniformization Theorem of recursion theory yields the existence of an $n${\hyp}ary partial recursive function $\varphi$ such that, whenever $x\in\mathbb{N}^n$ and the condition~(\ref{l}) holds for some $y\in\mathbb{N}$, then $x\in\mathrm{dom}(\varphi)$ and (\ref{l}) is satisfied by $y=\varphi(x)$. The equivalence of (\ref{l}) to the equality $f(x)=y$ for any $x\in\mathrm{dom}(f)$ and any $y\in\mathbb{N}$ implies that $f$ is a restriction of $\varphi$.

Suppose now the collection $\mathcal{A}$ is strongly join permitting for $\mathbf{C}$. Let $\mathcal{K}$ be an arbitrary subcollection of $\mathcal{A}$. We choose two distinct natural numbers $c_1$ and $c_2$ and consider the corresponding function $f$ defined as in Lemma~\ref{const}. The function $f$ belongs to $\mathbf{C}$, since $\mathrm{dom}(f)\subseteq\bigcup\mathcal{A}$ and $f\!\!\upharpoonright_{A}\in\mathbf{C}$ for any $A\in\mathcal{A}$.  Thus $f$ is a restriction of some $n${\hyp}ary partial recursive function~$\varphi$. Then $\varphi^{-1}(c_1)$ is a recursively enumerable subset of $\mathbb{N}^n$ containing $\bigcup\mathcal{K}\setminus\bigcup(\mathcal{A}\setminus\mathcal{K})$ and having an empty intersection with $\bigcup(\mathcal{A}\setminus\mathcal{K})\setminus\bigcup\mathcal{K}$.
\qed

\begin{remark}\label{proper}
{\em A recursively enumerable subset of $\mathbb{N}^n$ separating $\bigcup\mathcal{K}\setminus\bigcup(\mathcal{A}\setminus\mathcal{K})$ from $\bigcup(\mathcal{A}\setminus\mathcal{K})\setminus\bigcup\mathcal{K}$ obviously exists if $\mathcal{K}=\varnothing$ or $\mathcal{K}=\mathcal{A}$ (the empty set in the first case and $\mathbb{N}^n$ in the second one), therefore one may exclude these two cases in the condition formulated in Theorem~\ref{ppr}.}
\end{remark}

\begin{cor}\label{re}
Any finite collection of recursively enumerable subsets of $\mathbb{N}^n$ is strongly join permitting for $\mathbf{C}$.
\end{cor}

\proof
Let $\mathcal{A}$ be a finite collection of recursively enumerable subsets of $\mathbb{N}^n$, and $\mathcal{K}$ be a subcollection of $\mathcal{A}$. Then $\bigcup\mathcal{K}$ is a recursively enumerable set which contains $\bigcup\mathcal{K}\setminus\bigcup(\mathcal{A}\setminus\mathcal{K})$ and has an empty intersection with $\bigcup(\mathcal{A}\setminus\mathcal{K})\setminus\bigcup\mathcal{K}$.
\qed

\begin{remark}\label{simpler}
{\em In the situation from Corollary \ref{re}, the proof of Theorem \ref{ppr}, when carried out with $H_\mathcal{K}=\bigcup\mathcal{K}$, turns out to be unnecessarily complicated, because the condition \eqref{l} is equivalent to the much simpler condition
$$\exists A\in\mathcal{A}(x\in A\ \&\ \varphi_A(x)=y)$$
in this case.}
\end{remark}

\begin{cor}\label{core}
Let $\mathcal{A}$ be the collection $\{\mathbb{N}^n\setminus E\,|\,E\in\mathcal{E}\}$, where $\mathcal{E}$ is a finite collection of recursively enumerable subsets of $\mathbb{N}^n$. Then $\mathcal{A}$ is strongly join permitting for $\mathbf{C}$.
\end{cor}

\proof Let $\mathcal{K}$ be a subcollection of $\mathcal{A}$ distinct from $\mathcal{A}$, and let us set
$$H=\bigcap\{E\in\mathcal{E}\,|\,\mathbb{N}^n\setminus E\not\in\mathcal{K}\}.$$
Then $H$ is a recursively enumerable subset of $\mathbb{N}^n$ which contains $\bigcup\mathcal{K}\setminus\bigcup(\mathcal{A}\setminus\mathcal{K})$ and has an empty intersection with $\bigcup(\mathcal{A}\setminus\mathcal{K})\setminus\bigcup\mathcal{K}$.
\qed

\begin{remark}\label{rec} 
{\em Theorem \ref{ppr} remains true after skipping ``partial'' and replacing ``recursively enumerable'' with ``recursive'' in the definition of the class $\mathbf{C}$ (the potentially recursive $n${\hyp}ary functions are those ones which are restrictions of $n${\hyp}ary recursive functions, and clearly separability by means of a recursive set is the same thing as recursive separability). The main change in the proof consists in replacing the application of the Uniformization Theorem with an appropriate definition by cases. A similar version of the theorem concerning potential primitive recursiveness and primitive recursive separability holds too.}
\end{remark}

\section{Strongly join permitting collections \\for the class of the continuous partial functions \\between two given topological spaces}\label{S:continuous}

Let $\mathfrak{X}$ and $\mathfrak{Y}$ be topological spaces with carriers $X$ and $Y$, respectively. and $\mathbf{C}$ be the class of all continuous partial functions from $\mathfrak{X}$ to $\mathfrak{Y}$. 

\begin{thm}\label{tcont}
Let $\mathcal{A}$ be a finite collection of subsets of~$X$. If, for each subcollection $\mathcal{K}$ of~$\mathcal{A}$, some open set of $\mathfrak{X}$ separates $\bigcup\mathcal{K}\setminus\bigcup(\mathcal{A}\setminus\mathcal{K})$ from $\bigcup(\mathcal{A}\setminus\mathcal{K})\setminus\bigcup\mathcal{K}$ then the collection $\mathcal{A}$ is strongly join permitting for $\mathbf{C}$. In the case when there exists in $\mathfrak{Y}$ an open set different from $\varnothing$ and $Y$, the converse also holds.
\end{thm}

\proof
Suppose that, for any subcollection $\mathcal{K}$ of $\mathcal{A}$, some open set $H_\mathcal{K}$ of~$\mathfrak{X}$ separates \mbox{$\bigcup\mathcal{K}\setminus\bigcup(\mathcal{A}\setminus\mathcal{K})$} from $\bigcup(\mathcal{A}\setminus\mathcal{K})\setminus\bigcup\mathcal{K}$. We will prove that $\mathcal{A}$ is strongly join permitting for $\mathbf{C}$. Let $f$ be a function such that $\mathrm{dom}(f)\subseteq\bigcup\mathcal{A}$ and $f\!\upharpoonright_A\in\mathbf{C}$ for any $A\in\mathcal{A}$. Clearly, all values of $f$ belong to $Y$. We will show that $f\in\mathbf{C}$ by proving that, for any open set $V$ of $\mathfrak{Y}$, the set $f^{-1}(V)$ is the intersection of $\mathrm{dom}(f)$ with some open set of $\mathfrak{X}$. Let $V$ be an open set of $\mathfrak{Y}$. For any $A\in\mathcal{A}$ the set $(f\!\upharpoonright_A)^{-1}(V)$ is the intersection of $\mathrm{dom}(f\!\upharpoonright_A)$ with some open set $O_A$ of $\mathfrak{X}$. Let $x$ be an arbitrary element of $X$. We will show that $x\in f^{-1}(V)$ iff $x\in\mathrm{dom}(f)$ and
\begin{equation}\label{mn}
\exists\mathcal{K}\subseteq\mathcal{A}(x\in H_\mathcal{K}\ \&\ \forall A\in\mathcal{K}(x\in O_A)).
\end{equation}
Firstly, suppose that $x\in\mathrm{dom}(f)$ and $x\in H_\mathcal{K}\ \&\ \forall A\in\mathcal{K}(x\in O_A)$ for some subcollection~$\mathcal{K}$ of~$\mathcal{A}$. Since $x\in \bigcup\mathcal{A}$, statement (a) of Lemma \ref{lss} yields that $x\in\bigcup\mathcal{K}$, i.e.\@ $x\in A$ for some $A\in\mathcal{K}$. But if $x\in A$ then surely $x\in\mathrm{dom}(f\!\upharpoonright_{A})$. Since, in addition, $x\in O_A$, it follows that $x\in(f\!\upharpoonright_{A})^{-1}(V)$ and therefore $x\in f^{-1}(V)$. 
For the reasoning in the opposite direction, suppose that $x\in f^{-1}(V)$. Then, of course, $x\in\mathrm{dom}(f)$ and therefore $x\in \bigcup\mathcal{A}$. By statement~(b) of Lemma \ref{lss},  $x\in H_\mathcal{K}$, where $\mathcal{K}=\{A\in\mathcal{A}\,|\,x\in A\}$. Clearly, whenever $A$ belongs to this~$\mathcal{K}$, then $f\!\upharpoonright_{A}\!\!(x)=f(x)\in V$, consequently $x\in(f\!\upharpoonright_{A})^{-1}(V)$ and therefore $x\in O_A$.

It is clear now that $f^{-1}(V)$ is the intersection of $\mathrm{dom}(f)$ with the set of those $x\in X$ which satisfy the condition (\ref{mn}), and this set is an open set of $\mathfrak{X}$ thanks to the fact that $\mathcal{A}$ is a finite set.

Suppose now that the collection $\mathcal{A}$ is strongly join permitting for $\mathbf{C}$ and some open set of~$\mathfrak{Y}$ is different from $\varnothing$ and $Y$. Consider an arbitrary subcollection $\mathcal{K}$ of $\mathcal{A}$. Let $c_1$ and $c_2$ be elements of~$Y$ belonging to some open set $V$ of~$\mathfrak{Y}$ and to its complement, respectively. We consider the corresponding function $f$ defined as in Lemma~\ref{const}. The function $f$ belongs to~$\mathbf{C}$, since $\mathrm{dom}(f)\subseteq\bigcup\mathcal{A}$ and $f\!\!\upharpoonright_{A}\in\mathbf{C}$ for any $A\in\mathcal{A}$. Therefore $f^{-1}(V)$ is the intersection of $\mathrm{dom}(f)$ with some open set $H$ of $\mathfrak{X}$, and, since $f^{-1}(V)=f^{-1}(c_1)$, it is clear that $H$ contains $\bigcup\mathcal{K}\setminus\bigcup(\mathcal{A}\setminus\mathcal{K})$ and has an empty intersection with~$\bigcup(\mathcal{A}\setminus\mathcal{K})\setminus\bigcup\mathcal{K}$.\qed

\begin{remark}\label{propert}
{\em Since open sets of $\mathfrak{X}$ separating $\bigcup\mathcal{K}\setminus\bigcup(\mathcal{A}\setminus\mathcal{K})$ from $\bigcup(\mathcal{A}\setminus\mathcal{K})\setminus\bigcup\mathcal{K}$ obviously exist if $\mathcal{K}=\varnothing$ or $\mathcal{K}=\mathcal{A}$ (the empty set in the first case and the set $X$ in the second one), these two cases can be excluded in the condition formulated in Theorem~\ref{tcont}.}
\end{remark}

\begin{remark}\label{open}
{\em Analogs of Corollaries \ref{re} and \ref{core} hold and can be proved in a similar way, namely the statement that any finite collection of open sets of $\mathfrak{X}$ is strongly join permitting for~$\mathbf{C}$ and the same for any finite collection of closed sets of $\mathfrak{X}$. However, straightforward direct proofs of these two statements are well-known, and, moreover, the validity of the first of the statements is shown in such a direct way without using the finiteness assumption.}
\end{remark}

\section{Strongly join permitting collections \\for the class of the functions which are computable \\with regard to a given pair of sequences of sets}\label{S:computable}

Suppose $\mathcal{U}=\{U_k\}_{k\in\mathbb{N}}$ is a sequence of sets. Then we set $\mathcal{U}^{-1}(x)=\{k\in\mathbb{N}\,|\,x\in U_k\}$ for any $x\in\bigcup_{k=0}^\infty U_k$. The total enumerations of the set $\mathcal{U}^{-1}(x)$ will be called {\em$\mathcal{U}${\hyp}names} of $x$. The element $x$ will be called {\em$\mathcal{U}${\hyp}computable} if at least one of these enumerations is recursive; of course, the $\mathcal{U}${\hyp}computability of $x$ is equivalent to the recursive enumerability of the set~$\mathcal{U}^{-1}(x)$.

Let $D_0,D_1,D_2,\ldots$ be the canonical enumeration of the set of all finite subsets of~$\mathbb{N}$. We define the sequence $\hat{\mathcal{U}}=\{\hat{U}_m\}_{m\in\mathbb{N}}$ as follows: $\hat{U}_m=\bigcap_{k\in D_m}\!\!U_k$, assuming that $\bigcap_{k\in\varnothing}U_k$ equals $\bigcup_{k=0}^\infty U_k$. The collection $\hat{\mathcal{U}}$ is obviously closed under finite intersection. One easily checks that
\begin{equation}\label{hat}
x\in\hat{U}_m\Longleftrightarrow D_m\subseteq\mathcal{U}^{-1}(x)
\end{equation}
for any $m$ in $\mathbb{N}$ and any $x\in\bigcup_{k=0}^\infty U_k$. 

In the sequel, sequences $\mathcal{U}=\{U_k\}_{k\in\mathbb{N}}$ and $\mathcal{V}=\{V_l\}_{l\in\mathbb{N}}$ of sets are supposed to be given. A partial function $f$ from $\bigcup_{k=0}^\infty U_k$ to $\bigcup_{l=0}^\infty V_l$ will be called {\em$(\mathcal{U},\mathcal{V})${\hyp}computable} if a recursive operator $\Gamma$ exists such that, for any $x\in\mathrm{dom}(f)$, $\Gamma$ transforms all $\mathcal{U}${\hyp}names of $x$ into $\mathcal{V}${\hyp}names of $f(x)$. One proves that $f$ is $(\mathcal{U},\mathcal{V})${\hyp}computable iff an enumeration operator~$F$ exists such that $\mathcal{V}^{-1}(f(x))=F(\mathcal{U}^{-1}(x))$ for any $x\in\mathrm{dom}(f)$.\footnote{The proof can be done by carrying out certain reasonings from \cite[\S 9.7]{trfec} in an appropriate more formal way.} Thus $f$ is $(\mathcal{U},\mathcal{V})${\hyp}computable iff a recursively enumerable subset $W$ of $\mathbb{N}^2$ with the following property exists:
\begin{equation}\label{eo}
\forall x\in\mathrm{dom}(f)\left(\,\mathcal{V}^{-1}(f(x))=\{l\,|\,\exists m((m,l)\in W\ \&\ D_m\subseteq\mathcal{U}^{-1}(x))\}\,\right).
\end{equation}
By \eqref{hat}, the above property is equivalent to the following one:
$$\forall x\in\mathrm{dom}(f)\left(\,\mathcal{V}^{-1}(f(x))=\{l\,|\,\exists m((m,l)\in W\ \&\ x\in\hat{U}_m\}\,\right).$$
By rewriting it in the form
\begin{equation}\label{eoha}
\forall x\in\mathrm{dom}(f)\,\forall l\in\mathbb{N}\left(f(x)\in V_l\Leftrightarrow\exists m\left((m,l)\in W\ \&\ x\in\hat{U}_m\right)\,\right),
\end{equation}
we see that the considered property of $W$ is equivalent to being a $(\hat{\mathcal{U}},\mathcal{V})${\hyp}approximation system for $f$ in the sense of \cite[Definition 2.1]{ed} in the particular case when $\mathcal{U}$ and $\mathcal{V}$ are bases of some topological spaces (then obviously $\hat{\mathcal{U}}$ is also a base of the first of these spaces, and if they are $T_0$ spaces then the computability notions defined above coincide with the usual ones from the theory of computability in topological spaces).

We note that, whenever $c$ is a $\mathcal{V}${\hyp}computable element of $\bigcup_{l=0}^\infty V_l$, all partial functions from $\bigcup_{k=0}^\infty U_k$ to~$\{c\}$ are $(\mathcal{U},\mathcal{V})${\hyp}computable (since, for instance, the set $W=\mathbb{N}\times\mathcal{V}^{-1}(c)$ has the property~\eqref{eo} for any such function $f$).

The sets of the form $\bigcup_{m\in S}\hat{U}_m$, where $S$ is some recursively enumerable subset of $\mathbb{N}$, will be called {\em effective $\hat{\mathcal{U}}${\hyp}unions}. The union and the intersection of any two effective $\hat{\mathcal{U}}${\hyp}unions can be shown to be effective $\hat{\mathcal{U}}${\hyp}unions again.

\begin{lem}\label{pre-image}
Let $f$ be a $(\mathcal{U},\mathcal{V})${\hyp}computable partial functions from $\bigcup_{k=0}^\infty U_k$ to $\bigcup_{l=0}^\infty V_l$. Then, for any $l\in\mathbb{N}$, the set $f^{-1}(V_l)$ is the intersection of $\mathrm{dom}(f)$ with some effective $\hat{\mathcal{U}}${\hyp}union.
\end{lem}

\proof
By the $(\mathcal{U},\mathcal{V})${\hyp}computability of $f$, a recursively enumerable subset $W$ of $\mathbb{N}^2$ with the property \eqref{eoha} exists. This property is equivalent to the following one:
$$\forall l\in\mathbb{N}\left(f^{-1}(V_l)=\mathrm{dom}(f)\cap\left\{x\,\left|\,\exists m\left((m,l)\in W\ \&\ x\in\hat{U}_m\right)\right.\right\}\right).$$
Let $l$ be an arbitrary natural number, and let $S=\{m\,|\,(m,l)\in W\}$. Then $S$ is a recursively enumerable set of natural numbers, and the equality $f^{-1}(V_l)=\mathrm{dom}(f)\cap\bigcup_{m\in S}\hat{U}_m$ holds.
\qed

Let $\mathbf{C}$ be the class of all $(\mathcal{U},\mathcal{V})${\hyp}computable partial functions from $\bigcup_{k=0}^\infty U_k$ to $\bigcup_{l=0}^\infty V_l$.

\begin{thm}\label{tcomp}
Let $\mathcal{A}$ be a finite collection of subsets of~$\bigcup_{k=0}^\infty U_k$. If, for each subcollection $\mathcal{K}$ of $\mathcal{A}$, some effective $\hat{\mathcal{U}}${\hyp}union separates $\bigcup\mathcal{K}\setminus\bigcup(\mathcal{A}\setminus\mathcal{K})$ from $\bigcup(\mathcal{A}\setminus\mathcal{K})\setminus\bigcup\mathcal{K}$ then the collection $\mathcal{A}$ is strongly join permitting for~$\mathbf{C}$. If natural numbers $l_1$ and $l_2$ exist such that each of the sets $V_{l_1}$ and $V_{l_2}\setminus V_{l_1}$ contains some $\mathcal{V}${\hyp}computable element then the converse also holds.
\end{thm}

\proof
Suppose that, for any subcollection $\mathcal{K}$ of $\mathcal{A}$, some effective $\hat{\mathcal{U}}${\hyp}union $H_\mathcal{K}$ separates $\bigcup\mathcal{K}\setminus\bigcup(\mathcal{A}\setminus\mathcal{K})$ from $\bigcup(\mathcal{A}\setminus\mathcal{K})\setminus\bigcup\mathcal{K}$. We will prove that $\mathcal{A}$ is strongly join permitting for~$\mathbf{C}$. Let $f$ be a function such that $\mathrm{dom}(f)\subseteq\bigcup\mathcal{A}$ and $f\!\upharpoonright_A\in\mathbf{C}$ for any $A\in\mathcal{A}$. Clearly, $f$ is a partial function from $\bigcup_{k=0}^\infty U_k$ to $\bigcup_{l=0}^\infty V_l$. For any $A\in\mathcal{A}$, a recursively enumerable subset $W_A$ of~$\mathbb{N}^2$ can be chosen such that
\begin{equation}\label{mm}
l\in\mathcal{V}^{-1}(f(x))\Longleftrightarrow\exists m((m,l)\in W_A\,\&\,D_m\subseteq\mathcal{U}^{-1}(x)).
\end{equation}
for all $l\in\mathbb{N}$ and all $x\in A\cap\,\mathrm{dom}(f)$.
We will show that, for any $x\in\mathrm{dom}(f)$ and any $l\in\mathbb{N}$ the condition $l\in\mathcal{V}^{-1}(f(x))$ is equivalent to the condition
\begin{equation}\label{h}
\exists\mathcal{K}\subseteq\mathcal{A}\left(x\in H_\mathcal{K}\,\&\,\forall A\in\mathcal{K}\,\exists m\left((m,l)\in W_A\,\&\,D_m\subseteq\mathcal{U}^{-1}(x)\right)\right).
\end{equation}
Firstly, suppose that $x\in\mathrm{dom}(f)$ (hence $x\in \bigcup\mathcal{A}$), $l\in\mathbb{N}$, and $\mathcal{K}$ is a subcollection of $\mathcal{A}$ such that
$$x\in H_\mathcal{K}\,\&\,\forall A\in\mathcal{K}\,\exists m\left((m,l)\in W_A\,\&\,D_m\subseteq\mathcal{U}^{-1}(x)\right).$$
The statement~(a) of Lemma \ref{lss} implies that $x\in A$ for some $A\in\mathcal{K}$. Hence (by the equivalence \eqref{mm}) $l\in\mathcal{V}^{-1}(f(x))$.  For the reasoning in the opposite direction, suppose that \mbox{$l\in\mathcal{V}^{-1}(f(x))$.} By statement (b) of Lemma \ref{lss}, $x\in H_\mathcal{K}$, where $\mathcal{K}=\{A\in\mathcal{A}\,|\,x\in A\}$. On the other hand, again using the equivalence \eqref{mm}, we see that \mbox{$\exists m((m,l)\in W_A\,\&\,D_m\subseteq\mathcal{U}^{-1}(x))$} holds for any $\mathcal{A}$ in $\mathcal{K}$. 

For any subcollection $\mathcal{K}$ of $\mathcal{A}$, since $H_\mathcal{K}$ is an effective $\hat{\mathcal{U}}${\hyp}union, some recursively enumerable subset $S_\mathcal{K}$ of $\mathbb{N}$ exists such that $H_\mathcal{K}=\bigcup_{m\in S_\mathcal{K}}\!\hat{U}_m$. Making use also of the equivalence~\eqref{hat}, we may write the condition \eqref{h} in the form
$$\exists\mathcal{K}\subseteq\mathcal{A}\left(\exists m\in S_\mathcal{K}\left(D_m\subseteq\mathcal{U}^{-1}(x)\right)\,\&\,\forall A\in\mathcal{K}\,\exists m\left((m,l)\in W_A\,\&\,D_m\subseteq\mathcal{U}^{-1}(x)\right)\right).$$
The above condition is equivalent to
$$\exists \tilde{m}((\tilde{m},l)\in W\ \&\ D_{\tilde{m}}\subseteq\mathcal{U}^{-1}(x)),$$
where $W$ is the set of all $(\tilde{m},l)\in\mathbb{N}^2$ such that
$$\exists\mathcal{K}\subseteq\mathcal{A}\left(\exists m\in S_\mathcal{K}\left(D_m\subseteq D_{\tilde{m}}\right)\,\&\,\forall A\in\mathcal{K}\,\exists m\left((m,l)\in W_A\,\&\,D_m\subseteq D_{\tilde{m}}\right)\right).$$
Thus
$$\mathcal{V}^{-1}(f(x))=\{l\,|\,\exists \tilde{m}((\tilde{m},l)\in W\ \&\ D_{\tilde{m}}\subseteq\mathcal{U}^{-1}(x))\}$$
for all $x\in\mathrm{dom}(f)$. As it is easy to see, the set $W$ is recursively enumerable. Hence the function $f$ is $(\mathcal{U},\mathcal{V})${\hyp}computable.

Suppose now the collection $\mathcal{A}$ is strongly join permitting for $\mathbf{C}$, and $l_1,l_2$ are natural numbers such that $V_{l_1}$ and $ V_{l_2}\setminus V_{l_1}$ contain $\mathcal{V}${\hyp}computable elements $c_1$ and $c_2$, respectively. For an arbitrary subcollection $\mathcal{K}$ of $\mathcal{A}$, we consider the corresponding function $f$ defined as in Lemma~\ref{const}. By Lemma~\ref{const}, $\mathrm{dom}(f)\subseteq\bigcup\mathcal{A}$, and $f\!\!\upharpoonright_{A}$ is a constant function for any $A\in\mathcal{A}$. Taking into consideration the $\mathcal{V}${\hyp}computability of $c_1$ and $c_2$, we see that $f\!\!\upharpoonright_{A}\in\mathbf{C}$ for any \mbox{$A\in\mathcal{A}$.} Consequently $f\in\mathbf{C}$. By Lemma \ref{pre-image}, the set $f^{-1}(V_{l_1})$ is the intersection of $\mathrm{dom}(f)$ with some effective $\hat{\mathcal{U}}${\hyp}union $H$, and, since $f^{-1}(V_{l_1})=f^{-1}(c_1)$, it is clear that $H$ contains \mbox{$\bigcup\mathcal{K}\setminus\bigcup(\mathcal{A}\setminus\mathcal{K})$} and has an empty intersection with $\bigcup(\mathcal{A}\setminus\mathcal{K})\setminus\bigcup\mathcal{K}$.
\qed

\begin{remark}\label{rel}
{\em If $\mathcal{A}$ is a finite collection of sets relevant for $\mathbf{C}$ then the assumption in the first sentence of Theorem \ref{tcomp} is surely satisfied, because any set relevant for~$\mathbf{C}$ is a subset of $\bigcup_{k=0}^\infty U_k$ (the converse may turn out to be false even in the case when $\mathcal{U}$ and $\mathcal{V}$ are bases of $T_0$ spaces).}
\end{remark}

\begin{remark}\label{properc}
{\em Since effective $\hat{\mathcal{U}}${\hyp}unions separating $\bigcup\mathcal{K}\setminus\bigcup(\mathcal{A}\setminus\mathcal{K})$ from $\bigcup(\mathcal{A}\setminus\mathcal{K})\setminus\bigcup\mathcal{K}$ obviously exist if $\mathcal{K}=\varnothing$ or $\mathcal{K}=\mathcal{A}$ (namely $\varnothing$ and $\bigcup_{k=0}^\infty U_k$, respectively), one may exclude these two cases in the condition formulated in Theorem~\ref{tcomp}.}
\end{remark}

\begin{cor}\label{c}
Any finite collection of effective $\hat{\mathcal{U}}${\hyp}unions is strongly join permitting for~$\mathbf{C}$.
\end{cor}

\proof
If $\mathcal{A}$ is a finite collection of $\hat{\mathcal{U}}${\hyp}unions then, for each subcollection $\mathcal{K}$ of $\mathcal{A}$, the set $\bigcup\mathcal{K}$ is an effective $\hat{\mathcal{U}}${\hyp}union containing $\bigcup\mathcal{K}\setminus\bigcup(\mathcal{A}\setminus\mathcal{K})$ and having an empty intersection with $\bigcup(\mathcal{A}\setminus\mathcal{K})\setminus\bigcup\mathcal{K}$.
\qed

\begin{remark}\label{csimpler}
{\em In the situation from Corollary \ref{c} and its proof, the proof of Theorem \ref{tcomp}, when carried out with $H_\mathcal{K}=\bigcup\mathcal{K}$, turns out to be unnecessarily complicated, because the condition \eqref{h} is equivalent to the much simpler condition
$$\exists A\in\mathcal{A}\left(x\in A\ \&\ \exists m\left((m,l)\in W_A\,\&\,D_m\subseteq\mathcal{U}^{-1}(x)\right)\right)$$
in this case.}
\end{remark}

\begin{cor}\label{ccomp}
Let $\mathcal{A}$ be the collection $\{\bigcup_{k=0}^\infty U_k\setminus E\,|\,E\in\mathcal{E}\}$, where $\mathcal{E}$ is a finite collection of effective $\hat{\mathcal{U}}${\hyp}unions. Then $\mathcal{A}$ is strongly join permitting for $\mathbf{C}$.
\end{cor}

\proof
Let $\mathcal{K}$ be a subcollection of $\mathcal{A}$ distinct from $\mathcal{A}$, and let us set 
$$H=\bigcap\{E\in\mathcal{E}\,|\,\bigcup_{k=0}^\infty U_k\setminus E\not\in\mathcal{K}\}.$$
Then $H$ is an effective $\hat{\mathcal{U}}${\hyp}union which contains $\bigcup\mathcal{K}\setminus\bigcup(\mathcal{A}\setminus\mathcal{K})$ and has an empty intersection with $\bigcup(\mathcal{A}\setminus\mathcal{K})\setminus\bigcup\mathcal{K}$.
\qed

\begin{example}\label{arctan}
{\em The computability proof for the $\arctan$ function indicated in the introduction is not satisfactory enough from the point of view of numerical calculations due to the poor convergency rate of the series in the right-hand side of the equality \eqref{ml} for the values of $x$ in the interval $[-1,1]$ which are its endpoints or are near to them. The following proof based on Corollary \ref{c} is better from this point of view. Let $\mathcal{U}=\mathcal{V}$ be a computable enumeration of the set of all open intervals with rational endpoints in $\mathbb{R}$. Then $\mathbf{C}$ is the class of all computable unary real functions, and if $a$ and $b$ are rational numbers such that $0<b<a$ then the collection of the intervals $(-a,a)$, $(b,+\infty)$ and $(-\infty,-b)$ is strongly join permitting for $\mathbf{C}$. If, additionally, $a<2$ and $ab\ge 1$ then this can be used to prove the computability of the function $f(x)=\arctan x$ by proving the computability of its restrictions to these intervals in the following way. We prove the computability of $f\!\upharpoonright(-a,a)$ by using the equality
$$\arctan x=\arctan\frac{x}{2}+\arctan\frac{x}{2+x^2}$$
and applying \eqref{ml} to the two terms in its right-hand side, then we solve the problem for the other two intervals as we did in the introduction for the intervals $[1,+\infty)$ and $(-\infty,-1]$.} 
\end{example}

\begin{example}\label{le1}
{\em The statement in Example \ref{intervals} directly follows from Corollary \ref{ccomp} by applying it for the case of $\mathcal{U}$ and $\mathcal{V}$ as in Example \ref{arctan}.}
\end{example}

\begin{example}\label{le2}
{\em Let $n$ be a natural number, $\mathcal{U}$ be a computable enumeration of the set of all parallelotopes which are Cartesian products of $n+2$ open intervals with rational endpoints in~$\mathbb{R}$, and $\mathcal{V}$ be the same as in Example \ref{arctan}. Then $\mathbf{C}$ is the class of all computable $(n+2)${\hyp}argument real functions, and, by Corollary~\ref{ccomp}, the pair of sets $\left\{(x,t,y_1,\ldots,y_n)\in\mathbb{R}^{n+2}\,\big|\,x\le t\right\}$ and $\left\{(x,t,y_1,\ldots,y_n)\in\mathbb{R}^{n+2}\,\big|\,x\ge t\right\}$ is strongly join permitting for $\mathbf{C}$. In the case of $n=2$, this implies, for instance, the computability of the function
$$\mathrm{cases}(x,t,y,z)=\left\{\begin{array}{ll}y&\textnormal{if $x<t$ or $y=z$,}\\z&\textnormal{if $x>t$ or $y=z$}\end{array}\right.$$
considered in \cite[Subsection 1.1]{e}.}
\end{example}

\begin{remark}\label{2ppr}
{\em The results in this section cover as particular instances the ones about strongly join permitting collections for the class of the potentially partial recursive $n${\hyp}ary functions. These functions are exactly the $(\mathcal{U},\mathcal{V})${\hyp}computable ones, the recursively enumerable subsets of $\mathbb{N}^n$ are exactly the effective $\hat{\mathcal{U}}${\hyp}unions, and any natural number is $\mathcal{V}${\hyp}computable if $U_0,U_1,U_2,\ldots$ is an injective computable enumeration of the set of the one-element subsets of $\mathbb{N}^n$, and $V_l=\{l\}$ for any $l\in\mathbb{N}$.}
\end{remark}

\appendix
\section{Strongly join permitting families of sets}

Let us call a family of sets $\{A_i\}_{i\in I}$ strongly join permitting for a class $\mathbf{C}$ of functions if the corresponding collection of sets $\{A_i\,|\,i\in I\}$ is strongly join permitting for $\mathbf{C}$. The results we proved can be easily transferred from collections to families of sets (of course, a family of sets $\{A_i\}_{i\in I}$ is called finite if the index set $I$ is finite). .For instance, Theorem \ref{ppr} goes into the following statement, where $n$ is some positive integer, $\mathbf{C}$ is the class of all potentially partial recursive $n${\hyp}ary functions, and, for any subset $K$ of $I$, $A_K^\cup$ denotes the set $\bigcup_{i\in K}A_i$.

\begin{thm}
Let $\{A_i\}_{i\in I}$ be a finite family of subsets of $\mathbb{N}^n$. The family $\{A_i\}_{i\in I}$ is strongly join permitting for~$C$ iff, for each subset $K$ of $I$, some recursively enumerable subset of $\mathbb{N}^n$ separates \mbox{$A_K^\cup\setminus A_{I\setminus K}^\cup$} from $A_{I\setminus K}^\cup\setminus A_K^\cup$.
\end{thm}

\proof
Let $\mathcal{A}=\{A_i\,|\,i\in I\}$. Suppose that, for each subset $K$ of $I$, some recursively enumerable subset of $\mathbb{N}^n$ separates $A_K^\cup\setminus A_{I\setminus K}^\cup$ from $A_{I\setminus K}^\cup\setminus A_K^\cup$. Let $\mathcal{K}$ be an arbitrary subcollection of the collection $\mathcal{A}$, and let $K=\{i\in I\,|\,A_i\in\mathcal{K}\}$. Then $A_K^\cup=\bigcup\mathcal{K}$ and $A_{I\setminus K}^\cup=\bigcup(\mathcal{A}\setminus\mathcal{K})$, hence some recursively enumerable subset of $\mathbb{N}^n$ separates $\bigcup\mathcal{K}\setminus\bigcup(\mathcal{A}\setminus\mathcal{K})$ from $\bigcup(\mathcal{A}\setminus\mathcal{K})\setminus\bigcup\mathcal{K}$. Therefore, by Theorem \ref{ppr}, the collection $\mathcal{A}$ is strongly join permitting for $\mathbf{C}$, i.e.\@ $\{A_i\}_{i\in I}$ is strongly join permitting for~$C$. For the reasoning in the opposite direction, suppose that $\{A_i\}_{i\in I}$ is strongly join permitting for~$C$, i.e.\@ $\mathcal{A}$ is strongly join permitting for $\mathbf{C}$. Let $K$ be an arbitrary subset of $I$, and let $\mathcal{K}=\{A_i\,|\,i\in K\}$. Then $\mathcal{K}$ is a subcollection of $\mathcal{A}$, hence, by Theorem \ref{ppr}, some recursively enumerable subset $H$ of $\mathbb{N}^n$ separates $\bigcup\mathcal{K}\setminus\bigcup(\mathcal{A}\setminus\mathcal{K})$ from $\bigcup(\mathcal{A}\setminus\mathcal{K})\setminus\bigcup\mathcal{K}$. Since $A_K^\cup=\bigcup\mathcal{K}$ and $A_{I\setminus K}^\cup\supseteq\bigcup(\mathcal{A}\setminus\mathcal{K})$, the set $H$ separates also $A_K^\cup\setminus A_{I\setminus K}^\cup$ from $A_{I\setminus K}^\cup\setminus A_K^\cup$.
\qed

\section*{Acknowledgments}

Special thanks are due to an anonymous referee for raising a pertinent question concerning Theorem \ref{tcont} and to Georgi Dimov for a useful discussion on this subject.

\end{document}